\input amstex\input epsf\documentstyle{amsppt}\nologo\footline={}

\def\SU{\mathop{\text{\rm SU}}}
\def\tr{\mathop{\text{\rm tr}}}
\def\Im{\mathop{\text{\rm Im}}}
\def\End{\mathop{\text{\rm End}}}
\def\area{\mathop{\text{\rm area}}}
\def\G{\mathop{\text{\rm G}}}

\hsize450pt\vsize599pt\topmatter\title On `A characterization of
$\Bbb R$-Fuchsian groups acting on the complex hyperbolic
plane'\endtitle\author Sasha Anan$'$in\endauthor\subjclass57M50
(57S25)\endsubjclass\abstract We indicate a $\Bbb C$-Fuchsian
counter-example to the result with the above title announced at
http://www.maths.dur.ac.uk/events/Meetings/LMS/2011/GAL11/program.pdf
and prove a stronger statement.\endabstract\endtopmatter\document

\rightheadtext{}\leftheadtext{}

\centerline{\bf1.~Introduction}

\medskip

The following result

\smallskip

\noindent
`{\sl We prove that a complex hyperbolic non-elementary Kleinian
group\/ $G$ acting on two-dimensional complex hyperbolic space\/
$\bold H_\Bbb C^2$ is\/ $\Bbb R$-Fuchsian, that is, $G$ leaves
invariant a totally real plane in\/ $\bold H_\Bbb C^2$, if and only if
every loxodromic element of\/ $G$ is either hyperbolic or loxodromic
whose elliptic part is of order\/ $2$.}'

\smallskip

\noindent
is announced at
http://www.maths.dur.ac.uk/events/Meetings/LMS/2011/GAL11/program.pdf
as a content of a one-hour talk.

The assertion as it stands is wrong (see a $\Bbb C$-Fuchsian
counter-example in Section 3). The following theorem directly implies a
corrected statement.

\medskip

{\bf Theorem.} {\sl Let\/ $V$ be a\/ $\Bbb C$-linear space equipped
with a hermitian form\/ $\langle-,-\rangle$ of signature\/ $++-$ and
let\/ $G\le\SU V$ be a subgroup such that the trace\/ $\tr g$ of every
loxodromic element\/ $g\in G$ belongs to\/ $\Bbb R\delta_g$, where\/
$\delta_g^3=1$. Suppose that\/ $G$ contains a loxodromic element. Then
either there exists a\/ $1$-dimensional\/ $G$-stable\/
$\Bbb C$-subspace in $V$ or there exists a totally real\/
$3$-dimensional\/ $G$-stable\/ $\Bbb R$-subspace in\/ $V$.}

\bigskip

\centerline{\bf2.~Proof of Theorem}

\medskip

We assume that there is no $1$-dimensional $G$-stable $\Bbb C$-subspace
in $V$.

\medskip

{\bf2.1.}~First, suppose that $\tr G\subset\Bbb R$.

Let $W\le V$ be a $G$-stable $\Bbb R$-subspace in $V$. Then the
$\Bbb C$-subspaces~$\Bbb CW$, $W\cap iW$, and
$W^\perp:=\big\{v\in V\mid\langle v,W\rangle=0\big\}$ are obviously
$G$-stable. It follows that $\dim_\Bbb RW$ cannot equal

\smallskip

\noindent
$\bullet$ $1$ because, otherwise, $\Bbb CW$ is a $1$-dimensional
$G$-stable $\Bbb C$-subspace in $V$;

\noindent
$\bullet$ $2$ because, otherwise, $\dim_\Bbb C\Bbb CW$ equals $1$ or
$2$ and, in the latter case, $\dim_\Bbb CW^\perp=1$;

\noindent
$\bullet$ $4$ because, otherwise, either $W$ is a complex subspace with
$\dim_\Bbb CW^\perp=1$ or $W+iW=\Bbb CW=V$ and
$\dim_\Bbb R(W\cap iW)=2$, that is, $\dim_\Bbb C(W\cap iW)=1$;

\noindent
$\bullet$ $5$ because, otherwise, $W+iW=\Bbb CW=V$ and
$\dim_\Bbb R(W\cap iW)=4$, that is, $\dim_\Bbb C(W\cap iW)=2$.

\smallskip

Suppose that $\dim_\Bbb RW=3$. Let $g\in G$ be loxodromic. The
eigenvalues of $g$ are $1,r^{-1},r$, where $0,\pm1\ne r\in\Bbb R$.
Denote by $e_0,e_1,e_2\in V$ the corresponding eigenvectors, where
$e_0$ is positive and orthogonal to the isotropic $e_1,e_2$ such that
$c:=\langle e_1,e_2\rangle\ne0$. Since $W\cap iW=0$, there is no
$\Bbb C$-subspace in $W$. Therefore, $\dim_\Bbb R(W\cap\Bbb Ce_i)\le1$.
On the other hand, since the characteristic polynomial of $g$ equals
$(x-1)(x-r^{-1})(x-r)$, there is a basis of eigenvectors of $g$ in $W$.
Thus, we can assume that $e_0,e_1,e_2\in W$. Clearly, $W$ is totally
real if $c\in\Bbb R$. Suppose that $c\notin\Bbb R$. Then
$\Im\langle W,w\rangle=0$ for $w\in W$ is equivalent to
$w\in\Bbb Re_0$. For any $h\in G$, we have
$0=\Im\langle W,e_0\rangle=\Im\langle hW,he_0\rangle=\Im\langle
W,he_0\rangle$.
So,~$Ge_0\subset\Bbb Re_0$. A contradiction.

Suppose that $V$ has no proper $G$-stable $\Bbb R$-subspaces. Let
$A:=\Bbb RG$ denote the real span of $G$ and $D:=\End_AV$ denote the
division $\Bbb R$-algebra of endomorphisms of the simple $A$-module
$V$ (Schur's lemma). By Artin-Wedderburn theorem, a quotient algebra
of $A$ is isomorphic to $\End V_D$. Since $\dim_\Bbb RV=6$, we~have
$D=\Bbb R$ and $\dim_\Bbb R\End V_D=36$ or $D=\Bbb C$ and
$\dim_\Bbb C\End V_D=9$. On the other hand, $A\le\End_\Bbb CV$ and
$\dim_\Bbb C\End_\Bbb CV=9$. Hence, $A=\End_\Bbb CV$, which contradicts
$\tr A\subset\Bbb R$.

\medskip

{\bf2.2.}~Without loss of generality, we can assume that $G$ contains
a nontrivial cubic root of unity. Then there exists a loxodromic
$g\in G$ with $\tr g\in\Bbb R$. In a suitable basis $e_0,e_1,e_2$
with the Gram matrix
$\left[\smallmatrix 1&0&0\\0&0&1\\0&1&0\endsmallmatrix\right]$,
such a $g$ has the form
$g:=\left[\smallmatrix1&0&0\\0&r^{-1}&0\\0&0&r\endsmallmatrix\right]$,
where $0,\pm1\ne r\in\Bbb R$.

\medskip

{\bf2.3.~Remark.} {\sl Let\/ $g\in G$ be loxodromic with\/
$\tr g\in\Bbb R$ and let\/ $e_0,e_1,e_2\in V$ be eigenvectors of\/ $g$
with the Gram matrix\/
$\left[\smallmatrix 1&0&0\\0&0&1\\0&1&0\endsmallmatrix\right]$. Then,
for every\/ $h\in G$, there exists a cubic root of unity\/ $\delta$
such that\/
$\langle he_0,e_0\rangle,\langle he_1,e_2\rangle,\langle
he_2,e_1\rangle,\tr(g^nh)\in\Bbb R\delta$
for all\/ $n\in\Bbb Z$.}

\medskip

{\bf Proof.} It is easy to see that
$\tr(g^nh)=\langle he_0,e_0\rangle+r^{-n}\langle
he_1,e_2\rangle+r^n\langle he_2,e_1\rangle$. If
$\langle he_1,e_2\rangle\ne0$ or $\langle he_2,e_1\rangle\ne0$, then
$g^nh$ is loxodromic for sufficiently large $|n|$. Therefore,
$\langle he_0,e_0\rangle,\langle he_1,e_2\rangle,\allowmathbreak\langle
he_2,e_1\rangle\in\Bbb R\delta$
for a suitable cubic root of unity $\delta$. If
$\langle he_1,e_2\rangle=\langle he_2,e_1\rangle=0$, then
$h=\left[\smallmatrix-\varepsilon^{-2}&0&0\\0&0&a\varepsilon\\
0&a^{-1}\varepsilon&0\endsmallmatrix\right]$
with $a>0$ and $|\varepsilon|=1$. Since
$h^2=\left[\smallmatrix\varepsilon^{-4}&0&0\\0&\varepsilon^2&0\\
0&0&\varepsilon^2\endsmallmatrix\right]$
and
$\langle h^2e_1,e_2\rangle=\langle
h^2e_2,e_1\rangle=\varepsilon^2\ne0$,
we obtain $\varepsilon^2\in\Bbb R\overline\delta$, where $\delta^3=1$.
Again, we get
$\langle he_0,e_0\rangle,\langle he_1,e_2\rangle,\langle
he_2,e_1\rangle\in\Bbb R\delta$
$_\blacksquare$

\medskip

{\bf2.4.~Lemma.} {\sl Let\/ $g,h\in G$ be loxodromic with\/
$\tr g,\tr h\in\Bbb R$. Then\/ $\tr(gh)\in\Bbb R$.}

\medskip

{\bf Proof.} In some bases $e_0,e_1,e_2$ and $f_0,f_1,f_2$ with Gram
matrix $\left[\smallmatrix 1&0&0\\0&0&1\\0&1&0\endsmallmatrix\right]$,
we respectively have
$g=\left[\smallmatrix1&0&0\\0&r^{-1}&0\\0&0&r\endsmallmatrix\right]$
and
$h=\left[\smallmatrix1&0&0\\0&s^{-1}&0\\0&0&s\endsmallmatrix\right]$,
where $0,\pm1\ne r,s\in\Bbb R$. Let $g_{ij}:=\langle e_i,f_j\rangle$.
Then $e_i=g_{i0}f_0+g_{i2}f_1+g_{i1}f_2$ for $i=0,1,2$. By Remark 2.3,
for every $n\in\Bbb Z$, there exists some cubic root of unity
$\delta_n$ such that
$\langle h^ne_0,e_0\rangle,\langle h^ne_1,e_2\rangle,\langle
h^ne_2,e_1\rangle\in\Bbb R\delta_n$.
Taking $\delta$ such that $\delta_n=\delta$ for infinitely many $n$'s,
from
$$\langle h^ne_0,e_0\rangle=g_{00}\overline
g_{00}+s^{-n}g_{02}\overline g_{01}+s^ng_{01}\overline g_{02},$$
$$\langle h^ne_1,e_2\rangle=g_{10}\overline
g_{20}+s^{-n}g_{12}\overline g_{21}+s^ng_{11}\overline g_{22},\qquad
\langle h^ne_2,e_1\rangle=g_{20}\overline g_{10}+s^{-n}g_{22}\overline
g_{11}+s^ng_{21}\overline g_{12},$$
we obtain
$$g_{00}\overline g_{00},g_{02}\overline g_{01},g_{01}\overline
g_{02},g_{10}\overline g_{20},g_{12}\overline g_{21},g_{11}\overline
g_{22},g_{20}\overline g_{10},g_{22}\overline g_{11},g_{21}\overline
g_{12}\in\Bbb R\delta.$$
If $g_{11}\overline g_{22}=0$, then $e_2\ne e_1=f_1\ne f_2$ or
$e_1\ne e_2=f_2\ne f_1$ (the equalities and inequalities are meant in
the projective sense). Hence, $g_{12}\overline g_{21}\ne0$. We conclude
that $\delta=1$
$_\blacksquare$

\medskip

{\bf2.5.~Lemma.} {\sl Let\/ $g,h_1,h_2\in G$ be such that\/ $g$ is
loxodromic and\/ $\tr g,\tr(g^nh_1),\tr(g^nh_2)\in\Bbb R$ for all\/
$n\in\Bbb Z$. Then\/ $\tr(g^nh_1^{-1}h_2)\in\Bbb R$ for all\/
$n\in\Bbb Z$.}

\medskip

{\bf Proof.} Using the symmetry between $h_1,h_2$ and replacing
$h_1,h_2$ by $g^kh_1,g^kh_2$, if necessary, we~can assume (as in the
proof of Remark 2.3) that $h_2$ is loxodromic unless both $h_1,h_2$
have the type
$\left[\smallmatrix\pm1&0&0\\0&0&a\varepsilon\\
0&a^{-1}\varepsilon&0\endsmallmatrix\right]$
in the basis related to $g$, where $a>0$ and $\varepsilon^2=\mp1$. In
this particular case, $h_1^{-1}h_2$ is diagonal with coefficients in
$\Bbb R\cup\Bbb Ri$. By Remark 2.3, for some cubic root of unity
$\delta$, we have $\tr(g^nh_1^{-1}h_2)\in\Bbb R\delta$ for all
$n\in\Bbb Z$. Therefore, the mentioned coefficients have to be real.

So, we assume that $h_2$ is loxodromic with $\tr h_2\in\Bbb R$. Suppose
that $\tr(g^nh_1^{-1}h_2)\in\Bbb R\delta$ for all $n\in\Bbb Z$, where
$\delta^3=1$ and $\delta\ne1$. For some $m\in\Bbb Z$, we have
$0\ne\tr(g^mh_1^{-1}h_2)\in\Bbb R\delta$ as, otherwise, we are done.
Hence, by Remark 2.3, $\tr(g^mh_1^{-1}h_2^n)\in\Bbb R\delta$ for all
$n\in\Bbb Z$. In particular, $\tr(g^mh_1^{-1})\in\Bbb R\delta$, which
implies $\tr(g^mh_1^{-1})=0$.

Suppose that $g^kh_1^{-1}$ is loxodromic for some $k\in\Bbb Z$. As
in the proof of Remark 2.3, we conclude that $g^nh_1^{-1}$ is
loxodromic for all sufficiently large/small $n$. By Lemma 2.4,
$\tr(g^nh_1^{-1}h_2)\in\Bbb R$ for all such $n$'s, implying
$\tr(g^nh_1^{-1}h_2)=0$, a contradiction.

So, $h_1^{-1}$ is of the type
$\left[\smallmatrix\pm1&0&0\\0&0&a\varepsilon\\
0&a^{-1}\varepsilon&0\endsmallmatrix\right]$.
This contradicts $\tr(g^mh_1^{-1})=0$
$_\blacksquare$

\medskip

{\bf2.6.}~By Lemma 2.5,
$H:=\big\{h\in G\mid\tr(g^nh)\in\Bbb R\text{ for all }n\in\Bbb Z\big\}$
is a subgroup in $G$. Obviously, $G$ is generated by $H$ and the cubic
roots of unity. It suffices to deal with $H$ in place of $G$. In other
words, we can assume that $\tr G\subset\Bbb R$.

\bigskip

\centerline{\bf3.~Counter-example}

\medskip

\vskip15pt

\noindent
$\vcenter
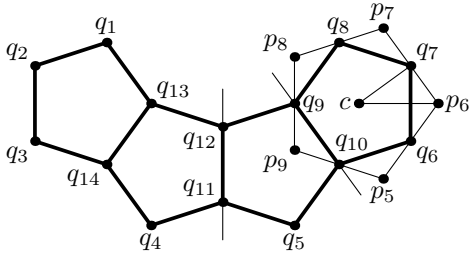
{\hbox{\epsfbox{Picts.3}}}$

\leftskip185pt

\vskip-108pt

Let $\Delta(c,p_6,q_7)$ be a geodesic triangle in the hyperbolic plane
with the corresponding interior angles $\frac\pi5,\frac\pi5,\frac\pi2$.
The area of this triangle equals $\frac\pi{10}$. Taking $10$ congruent
triangles with common vertex $c$, we obtain a pentagon with area
$\area(p_5,p_6,p_7,p_8,p_9)=\pi$. By [ABG], the reflections $R(q_i)$ in
the middle points $q_i$, $i=6,7,8,9,10$, of the sides of the pentagon
satisfy the relation $R(q_{10})R(q_9)R(q_8)R(q_7)R(q_6)=\pm1$ in
$\SU(1,1)$ and provide a discrete group $H_5$. Note that, by~the
definition from [ABG], we have
$R(q)x:=i\Big(x-2\displaystyle\frac{\langle x,q\rangle}{\langle
q,q\rangle}q\Big)$. 

\leftskip0pt

Denote $Q(q):=-iR(q)$ (in the complex hyperbolic plane,
$Q(q)\in\SU V$). We consider $3$ more copies of the pentagon
$P(q_6,q_7,q_8,q_9,q_{10})$, namely:
$P(q_5,q_{10},q_9,q_{12},q_{11})$,
$P(q_4,q_{11},q_{12},q_{13},q_{14})$, and
$P(q_1,q_2,q_3,q_{14},q_{13})$. The geodesics
$\G{\prec}q_9,q_{10}{\succ}$ and $\G{\prec}q_{11},q_{12}{\succ}$
are ultraparallel (this can be shown with the help of SEs; see [ABG]).
The geodesics of this type separate the four pentagons, so that we have
exactly what is drawn on the picture. Since
$Q(q_{10})Q(q_9)Q(q_8)Q(q_7)Q(q_6)=
\left[\smallmatrix-1&0&0\\0&
\pm(-i)^5&0\\0&0&\pm(-i)^5\endsmallmatrix\right]=\left[\smallmatrix-1
&0&0\\0&\mp i&0\\0&0&\mp i\endsmallmatrix\right]$
and $Q(q)Q(q)=1$ in $\SU V$, we have
$$1=\left[\smallmatrix-1&0&0\\0&\mp i&0\\0&0&\mp
i\endsmallmatrix\right]^4=
\big(Q(q_8)Q(q_7)Q(q_6)Q(q_{10})Q(q_9)\big)\cdot
\big(Q(q_9)Q(q_{10})Q(q_5)Q(q_{11})Q(q_{12})\big)\cdot$$
$$\cdot\big(Q(q_{12})Q(q_{11})Q(q_4)Q(q_{14})Q(q_{13})\big)\cdot
\big(Q(q_{13})Q(q_{14})Q(q_3)Q(q_2)Q(q_1)\big)=$$
$$=Q(q_8)Q(q_7)Q(q_6)Q(q_5)Q(q_4)Q(q_3)Q(q_2)Q(q_1).$$
By [ABG], we obtain a $\Bbb C$-Fuchsian faithful and discrete
representation of $H_8$ and, hence, a $\Bbb C$-Fuchsian faithful and
discrete representation of the fundamental group $G_8$ of a surface of
genus $3$. As $G_8$ consists of all words of even length in the
$Q(q_i)$'s, $i=1,2,3,4,5,6,7,8$, every element $I\in G_8$ has the form
$I=\left[\smallmatrix1&0&0\\0&r^{-1}\alpha&0\\
0&0&r\alpha\endsmallmatrix\right]$
in a suitable basis $e_0,e_1,e_2$, where $e_1,e_2$ are isotropic points
in the complex geodesic, $e_0$~is its polar point, $r>0$, and
$|\alpha|=1$. Since $I\in\SU V$, we obtain $\alpha=\pm1$.

\bigskip

\centerline{\bf4.~References}

\medskip

[ABG]  S.~Anan$'$in, E.~C.~Bento~Gon\c calves, {\it A hyperelliptic
view on Teichm\"uller space. {\rm I},} preprint
http://arxiv.org/abs/0709.1711

\enddocument